%% file: ellipsoids-CR-umbilics.tex
\documentclass[12pt,twoside,reqno,openany]{amsart}

\usepackage{amsbsy,amscd,amsfonts,amsmath,amssymb,amsthm,color,
fancybox,fancyhdr,footmisc,graphics,graphicx,ifthen,mathrsfs,
multicol,pdfpages,rotating,times,wasysym}

\usepackage[dvipsnames,svgnames,x11names]{xcolor}

\usepackage[all]{xy}
\usepackage[utf8]{inputenc}
\usepackage[T1]{fontenc}
\usepackage{mathtools}
\newtagform{EngelLie}[\scriptstyle]{$}{$}
\makeatletter\newcommand{\leqnomode}{\tagsleft@true}
\newcommand{\reqnomode}{\tagsleft@false}\makeatother

\input macros.tex


\begin{document}

$\:$

\bigskip\bigskip\bigskip\bigskip\bigskip

\begin{center}

{\large\bf
Parametric CR-umbilical Locus of Ellipsoids in $\C^2$
}
\label{CR-umbilics-ellipsoids}

\bigskip\bigskip

Wei-Guo {\sc Foo}, Jo\"el {\sc Merker}, The-Anh {\sc Ta}

\smallskip

{\large\footnotesize\sf Departement of Mathematics}

{\large\footnotesize\sf 
Orsay University}

{\large\footnotesize\sf 
Paris, France}

\end{center}\bigskip

\bigskip

%%%%%%%%%%%%%%%%%%%%%%%%%%%%%%%%%%%%%%%%%%%%%%%%%%%%%%%%%%%%%%%%%%%%%%
\Section{\bf Introduction}
\label{introduction}
\HEAD{\ref{introduction}.~{\sf Introduction}
}{
Wei-Guo {\sc Foo}, Jo\"el {\sc Merker}, The-Anh {\sc Ta}
Orsay University, Paris, France}
%%%%%%%%%%%%%%%%%%%%%%%%%%%%%%%%%%%%%%%%%%%%%%%%%%%%%%%%%%%%%%%%%%%%%%

\CITATION{
Different from the situation in the classical Differential Geometry,
except in the trivial spherical case, where $\mathcal{S}$ or 
$\mathcal{P} \equiv 0$, computing umbilical points seems to be a
very difficult problem. This is because the explicit formula for 
the fundamental Cartan-Chern-Moser curvature tensors is too 
complicated.
\hfill
Xiaojun {\sc Huang}, Shanyu {\sc Ji}, 
\cite{Huang-Ji-2007}
}

In 1932, \'Elie Cartan~\cite{Cartan-1932-I, Cartan-1932-II,
Cartan-1932-III} showed that
a local real-analytic ($\mathcal{C}^\omega$)
hypersurface $M^3 \subset \C^2$ is determined 
up to local biholomorphic equivalence by a single
invariant function:
\[
\ICartan^M
\colon\ \ \
M
\,\longrightarrow\,
\C,
\]
together with its (covariant) derivatives with respect to
a certain coframe of differential $1$-forms on an $8$-dimensional
principal bundle $P^8 \longrightarrow M$.
In coordinates $(z,w) = \big( x + \isqrt\, y,\,
u + \isqrt\, v\big)$ on $\C^2$, whenever $M$ is:

\smallskip\noindent$\bullet$\,
either a {\sl complex graph:}
\[
\big\{
(z,w)\in\C^2
\colon\,
w=\Theta\big(z,\overline{z},\overline{w}\big)
\big\},
\]

\smallskip\noindent$\bullet$\,
or a {\sl real graph:}
\[
\big\{
(z,w)\in\C^2
\colon\,
v=\varphi(x,y,u)
\big\},
\]

\smallskip\noindent$\bullet$\,
or represented in {\sl implicit form:}
\[
\big\{
(z,w)\in\C^2
\colon\,
\rho\big(z,w,\overline{z},\overline{w}\big)
=
0
\big\},
\]

\smallskip\noindent
it is known that $\ICartan^M$ depends on the respective $6$-jets:
\[
J_{z,\overline{z},\overline{w}}^6\,
\Theta,
\ \ \ \ \ \ \ \ \ \ \ \ \ \ \
J_{x,y,u}^6\,
\varphi,
\ \ \ \ \ \ \ \ \ \ \ \ \ \ \
J_{z,w,\overline{z},\overline{w}}^6\,\rho.
\]

The {\sl invariancy} of $\ICartan^M$ means that, for any local
biholomorphism $h \colon \C^2 \longrightarrow \C^2$, setting $M' :=
h(M)$, it holds at every point $p \in M$ that:
\[
\ICartan^{M'}
\big(h(p)\big)
\,=\,
\nu(p)\,
\ICartan^M(p)
\eqno
{\scriptstyle{(\forall\,p\,\in\,M)}},
\]
for some nowhere vanishing (local) function $\nu \colon M
\longrightarrow \C \backslash \{0\}$.  This guarantees that the locus
of {\sl CR-umbilical} points:
\[
\UmbCR(M)
\,:=\,
\big\{
p\in M
\colon\,
\ICartan^M(p)
=
0
\big\}
\]
is intrinsic. Furthermore, when $M$ is connected, it is well known
that $\UmbCR(M)$ contains an open set $\emptyset \neq V \subset M$ if
and only if $M$ is {\sl spherical}, in the sense of being locally
biholomorphic to the unit sphere $S^3 \subset \C^2$.

In 1974, Chern-Moser~{\cite{Chern-Moser-1974}} raised the problem
whether $\emptyset \neq \UmbCR(M)$ for compact Levi nondegenerate
$\mathcal{C}^\omega$ hypersurfaces $M^{2\NN-1} \subset \C^\NN$
when $\NN \geqslant 2$. 
This (simple!) paper attacks the more specific:

\begin{Question}
{\sl Can $\UmbCR(M)$ be described explicity?}
\end{Question}

But because $\ICartan^M$ is `{\sl too complicated}'
as confirmed in~{\cite{Merker-Sabzevari-2012, Merker-Sabzevari-2014}},
the question is nontrivial even in simplest nonspherical
examples like {\em e.g.} real ellipsoids introduced
and studied by Webster in~{\cite{Webster-1977,
Webster-2000}}.

In $\C^{\NN \geqslant 2} \cong \R^{2\NN \geqslant 4}$
equipped with coordinates $\zaux_i = \xaux_i + \isqrt\,
\yaux_i$, an {\sl 
ellipsoid} is the image of the unit sphere:
\[
S^{2\NN-1}
\,:=\,
\big\{
\zaux\in\C^\NN
\colon\,
\vert\zaux_1\vert^2
+\cdots+
\vert\zaux_\NN\vert^2
=
1
\big\},
\]
through a real affine transformation of $\R^{2\NN}$, hence 
has equation the form:
\leqnomode\usetagform{default}
\begin{align}
\sum_{1\leqslant i\leqslant\NN}\,
\big(
\alpha_i\,\xaux_i^2
+
\beta_i\,\yaux_i^2
\big)
\,=\,
1,
\tag{${\sf E}_{\alpha,\beta}$}
\end{align}
with real constants $\alpha_i \geqslant \beta_i > 0$\,\,---\,\,replace
$\zaux_i \longmapsto \isqrt\, \zaux_i$ if necessary.

The complex geometry of ellipsoids
(Segre varieties, dynamics) began in Webster's
seminal article~{\cite{Webster-1977}}, in which it was verified
that two ellipsoids ${\sf E}_{\alpha,\beta} \cong {\sf E}_{
\alpha', \beta'}$ are biholomorphically equivalent if and only
if up to permutation:
\[
\frac{\alpha_i-\beta_i}{\alpha_i+\beta_i}
\,=\,
\frac{\alpha_i'-\beta_i'}{\alpha_i'+\beta_i'}
\eqno
{\scriptstyle{(1\,\leqslant\,i\,\leqslant\,\NN)}}.
\]
Replacing $\zaux_i \longmapsto \frac{1}{\sqrt{\beta_i}}\, \zaux_i$
and setting $a_i := \frac{\alpha_i}{\beta_i}$, whence
$1 \leqslant a_i$, leads to a convenient representation:
\leqnomode\usetagform{default}
\begin{align}
\sum_{1\leqslant i\leqslant\NN}\,
\big(
a_i\,\xaux_i^2
+
\yaux_i^2
\big)
\,=\,
1.
\tag{${\sf E}_{a_1,\dots,a_\NN}$}
\end{align}

Yet an alternative view, due to Webster in~{\cite{Webster-2000}}, is:
\leqnomode\usetagform{default}
\begin{align}
\label{A-Webster-ellipsoid}
\sum_{1\leqslant i\leqslant\NN}\,
\big(
\zaux_i\overline{\zaux}_i
+
A_i\,
(\zaux_i^2+\overline{\zaux}_i^2)
\big)
\,=\,
1,
\tag{${\sf E}_{A_1,\dots,A_\NN}$}
\end{align}
obtained by setting $A_i := \frac{a_i-1}{2a_i+2}$, whence $0
\leqslant A_i < \frac{1}{2}$, so that $a_i = \frac{1+2A_i}{1-2A_i}$,
then by changing coordinates $\zaux_i =: \sqrt{1-2A_i}\,
\zaux_i'$, and then by dropping primes.

In $\C^\NN$ when $\NN \geqslant 3$, what corresponds
to the invariant $\ICartan^M$ is the Hachtroudi-Chern tensor
$S_{\rho\sigma}^{\alpha\beta}$ with indices $1 \leqslant \alpha,
\beta, \rho, \sigma \leqslant \NN$, and the concerned
CR-umbilical locus:
\[
\UmbCR(M)
\,:=\,
\big\{
p\in M
\colon\,
S_{\rho\sigma}^{\alpha\beta}(p)
=
0,\,\,
\forall\,\alpha,\beta,\rho,\sigma
\big\},
\]
is known, through local biholomorphisms $h \colon \C^\NN 
\longrightarrow \C^\NN$ as above, to enjoy
\[
h\big(
\UmbCR(M)
\big)
\,=\,
\UmbCR
\big(h(M)\big).
\]

\begin{Theorem}
{\rm ({\cite{Webster-2000}})} 
In $\C^{\NN \geqslant 3}$, if $0 < A_1 < \cdots < A_\NN < \frac{1}{2}$,
then:
\[
\emptyset
\,=\,
\UmbCR
\big(
{\sf E}_{A_1,\dots,A_\NN}
\big).
\eqno\qed
\]
\end{Theorem}

This motivated Huang-Ji in~{\cite{Huang-Ji-2007}} to study the
question for compact $\mathcal{C}^\omega$ hypersurfaces $M \subset
\C^2$. If $M = \{ \rho = 0\}$, the expected dimension of:
\[
\UmbCR(M)
\,=\,
\big\{
0
=
\rho
=
\Re\,\ICartan
=
\Im\,\ICartan
\big\}
\]
should be $4 - 3 = 1$, although this is not rigorous, for $\R$ is not
algebraically closed!

\begin{Theorem}
{\rm (Implicitly proved in~{\cite{Huang-Ji-2007}})}
Every real ellipsoid ${\sf E}_{a,b} \subset \C^2$ of equation:
\[
a\,x^2+y^2+b\,u^2+v^2
\,=\,
1
\eqno
{\scriptstyle{(a\,\geqslant\,1,\,b\,\geqslant\,1,\,
(a,b)\,\neq\,(1,1))}}
\]
enjoys:
\[
\dim_\R\,
\UmbCR(M)
\,\geqslant\,
1.
\eqno\qed
\]
\end{Theorem}

In other words, it contains at least some (real algebraic!) curve.

What curve? Simple? Complicated? 

Can what follows be considered as a satisfactory answer?

\begin{Theorem}
\label{a-b-cos-sin}
For every real numbers $a \geqslant 1$, $b \geqslant 1$
with $(a,b) \neq (1,1)$,  
the curve parametrized by $\theta \in \R$ valued in
$\C^2 \cong \R^4${\em :}
\[
\gamma
\colon
\ \ \
\theta
\,\,\,\longmapsto\,\,\,
\big(
x(\theta)+\isqrt\,y(\theta),\,\,
u(\theta)+\isqrt\,v(\theta)
\big)
\]
with components:
\[
\aligned
x(\theta)
&
\,:=\,
{\textstyle{\sqrt{\frac{a-1}{a\,(ab-1)}}}}\,
\cos\,\theta,
\ \ \ \ \ \ \ \ \ \ \ \ \ \ \ \ \ \ \ \ \ \ \ \ \ \ \ \ \ \ 
y(\theta)
\,:=\,
{\textstyle{\sqrt{\frac{b\,(a-1)}{ab-1}}}}\,
\sin\,\theta,
\\
u(\theta)
&
\,:=\,
{\textstyle{\sqrt{\frac{b-1}{b\,(ab-1)}}}}\,
\sin\,\theta,
\ \ \ \ \ \ \ \ \ \ \ \ \ \ \ \ \ \ \ \ \ \ \ \ \ \ \ \ \ \ 
v(\theta)
\,:=\,
-\,
{\textstyle{\sqrt{\frac{a\,(b-1)}{ab-1}}}}\,
\cos\,\theta,
\endaligned
\]
has image contained in the CR-umbilical locus:
\[
\gamma(\R)
\,\subset\,
\UmbCR
\big({\sf E}_{a,b}\big)
\,\subset\,
{\sf E}_{a,b}
\]
of the ellipsoid ${\sf E}_{a,b} \subset \C^2$
of equation $a\,x^2+y^2+b\,u^2+y^2 = 1$.
\end{Theorem}

In other words:
\[
\ICartan^{{\sf E}_{a,b}}\big(\gamma(\theta)\big)
\,=\,
0
\eqno
{\scriptstyle{(\forall\,\theta\,\in\,\R)}}.
\]

As is known for ellipsoids, Cartan's invariant
$\ICartan^{{\sf E}_{a,b}}$ exhibits a high complexity,
{\em e.g.} $\sim\,40\,000$ terms in~{\cite{Merker-Sabzevari-2012}}.
So this theorem might be interpreted as a somewhat
unexpectedly nice and simple description of
$\UmbCR \big({\sf E}_{a,b}\big)$!

All computations of this paper were done by hand.

%%%%%%%%%%%%%%%%%%%%%%%%%%%%%%%%%%%%%%%%%%%%%%%%%%%%%%%%%%%%%%%%%%%%%%
\Section{\bf Explicit Expression of Cartan's CR-Invariant
$\mathfrak{I}$}
\label{explicit-Cartan-invariant}
\HEAD{\ref{explicit-Cartan-invariant}.~{\sf Explicit Expression 
of Cartan's CR-Invariant $\mathfrak{I}$}
}{
Wei-Guo {\sc Foo}, Jo\"el {\sc Merker}, The-Anh {\sc Ta}
Orsay University, Paris, France}
%%%%%%%%%%%%%%%%%%%%%%%%%%%%%%%%%%%%%%%%%%%%%%%%%%%%%%%%%%%%%%%%%%%%%%

In $\C^2$ equipped with coordinates $(z,w) = \big( x + \isqrt\,y,\,
u + \isqrt\, v\big)$, consider a connected real-analytic
($\mathcal{C}^\omega$) $3$-dimensional hypersurface:
\[
M^3
\,:=\,
\big\{
(z,w)\in\C^2
\colon\,
\rho(z,w,\overline{z},\overline{w})
=
0
\big\},
\]
with $\overline{\rho} = \rho$, and with $d\rho \big\vert_M$
never zero. Local or global $M$, compact or open,
bounded or unbounded, can be equally treated.

The two vector fields:
\[
L
\,:=\,
-\,\rho_w\,
\frac{\partial}{\partial z}
+
\rho_z\,
\frac{\partial}{\partial w}
\ \ \ \ \ \ \ \ \ \ \ \ \
\text{\rm and}
\ \ \ \ \ \ \ \ \ \ \ \ \
\overline{L}
\,:=\,
-\,\rho_{\overline{w}}\,
\frac{\partial}{\partial\overline{z}}
+
\rho_{\overline{z}}\,
\frac{\partial}{\partial\overline{w}}
\]
generate $T^{1,0}M$ and $T^{0,1}M$.

If $h \colon \C^2 \longrightarrow \C^2$ is a local biholomorphism:
\[
(z,w)
\,\longmapsto\,
\big(f(z,w),g(z,w)\big)
\,=:\,
(z',w'),
\]
if $M = \{ \rho = 0\}$ and $M' = \{ \rho' = 0\}$ are
two $\mathcal{C}^\omega$ 
hypersurfaces, if $h(M) \subset M'$, there is a
nowhere vanishing function $\mu \colon M \longrightarrow
\C \backslash \{0\}$ such that:
\[
\mu(z,w,\overline{z},\overline{w})\,
\rho(z,w,\overline{z},\overline{w})
\,\,\equiv\,\,
\rho'
\big(
f(z,w),g(z,w),
\overline{f}(\overline{z},\overline{w}),
\overline{g}(\overline{z},\overline{w})
\big),
\]
whence in $\C\{z,w,\overline{z},\overline{w}\}$ (exercise):
\[
\mu\,
\Big(
-\,\rho_w\,
\frac{\partial}{\partial z}
+
\rho_z\,
\frac{\partial}{\partial w}
\Big)
\,\,=\,\,
\big(f_zg_w-f_wg_z\big)\,
\Big(
-\,\rho_{w'}'\,
\frac{\partial}{\partial z'}
+
\rho_{z'}'\,
\frac{\partial}{\partial w'}
\Big).
\]

Furthermore, the {\sl Levi determinant:}
\[
\aligned
\Levi(\rho)
\,:=\,
&\,
-\,
\left\vert
\begin{array}{ccc}
0 & \rho_z & \rho_w
\\
\rho_{\overline{z}} & \rho_{z\overline{z}} & \rho_{w\overline{z}}
\\
\rho_{\overline{w}} & \rho_{z\overline{w}} & \rho_{w\overline{w}}
\end{array}
\right\vert
\\
\,=\,
&\,
\rho_{\overline{z}}\rho_z\rho_{w\overline{w}}
-
\rho_{\overline{z}}\rho_w\rho_{z\overline{w}}
-
\rho_{\overline{w}}\rho_z\rho_{\overline{z}w}
+
\rho_{\overline{w}}\rho_w\rho_{z\overline{z}},
\endaligned
\]
enjoys (exercise):
\[
\mu^3\,
{\sf L}(\rho)
\,=\,
\big(f_zg_w-f_wg_z\big)\,
\big(\overline{f}_{\overline{z}}
\overline{g}_{\overline{w}}
-
\overline{f}_{\overline{w}}
\overline{g}_{\overline{z}}
\big)\,
{\sf L}(\rho')
\eqno
{\scriptstyle{(\text{\rm on}\,M)}}.
\]

\begin{Definition}
A smooth hypersurface $M^3 \subset \C^2$ is called 
{\sl Levi nondegenerate} at a point $p \in M$ if:
\[
0
\,\neq\,
{\sf L}(p).
\]
\end{Definition}

From now on, all $M$ will be assumed smooth and Levi nondegenerate 
at every point, without further mention.

When $0 \neq \rho_w(p) = \rho_{\overline{w}}(p)$ at a point 
$p = (z_p, w_p) \in M$, the implicit function theorem represents
$M$ as a complex graph:
\[
w
\,=\,
\Theta
\big(z,\overline{z},\overline{w}\big)
\ \ \ \ \ \ \ \ \ \ \ \ \
\text{\rm or equivalently:}
\ \ \ \ \ \ \ \ \ \ \ \ \
\overline{w}
\,=\,
\overline{\Theta}
\big(
\overline{z},z,w
\big),
\]
in terms of a $\mathcal{C}^\omega$ defining function $\Theta$.
A similar graphed representation exists at points $q = (z_q, w_q)
\in M$ at which $0 \neq \rho_z(q) = \rho_{\overline{z}}(q)$.

Differentiating the identity:
\[
0
\,\equiv\,
\rho
\big(
z,
\Theta\big(z,\overline{z},\overline{w}\big),
\overline{z},\overline{w}
\big)
\eqno
{\scriptstyle{(\text{\rm in}\,\C\{z,\overline{z},\overline{w}\})}},
\]
once with respect to $z$, $\overline{z}$, $\overline{w}$ yields:
\[
\aligned
0
&
\,\equiv\,
\rho_z
+
\Theta_z\,\rho_w,
\notag
\\
0
&
\,\equiv\,
\rho_{\overline{z}}
+
\Theta_{\overline{z}}\,
\rho_w,
\\
0
&
\,\equiv\,
\rho_{\overline{w}}
+
\Theta_{\overline{w}}\,
\rho_w,
\endaligned
\]
and next twice with respect to $zz$, $z\overline{z}$,
$z\overline{w}$, $\overline{z}\overline{z}$, 
$\overline{z}\overline{w}$, $\overline{w}\overline{w}$ gives:
\leqnomode\usetagform{default}
\begin{align}
\label{2-rho-Theta}
0
&
\,\equiv\,
\rho_{zz}
+
2\,\Theta_z\,
\rho_{zw}
+
\Theta_z\,\Theta_z\,
\rho_{ww}
+
\Theta_{zz}\,
\rho_w,
\notag
\\
0
&
\,\equiv\,
\rho_{z\overline{z}}
+
\Theta_z\,
\rho_{\overline{z}w}
+
\Theta_{\overline{z}}\,
\rho_{zw}
+
\Theta_z\,\Theta_{\overline{z}}\,
\rho_{ww}
+
\Theta_{z\overline{z}}\,
\rho_w,
\notag
\\
0
&
\,\equiv\,
\rho_{z\overline{w}}
+
\Theta_z\,
\rho_{w\overline{w}}
+
\Theta_{\overline{w}}\,
\rho_{zw}
+
\Theta_z\,\Theta_{\overline{w}}\,
\rho_{ww}
+
\Theta_{z\overline{w}}\,
\rho_w,
\\
0
&
\,\equiv\,
\rho_{\overline{z}\overline{z}}
+
2\,\Theta_{\overline{z}}\,
\rho_{\overline{z}w}
+
\Theta_{\overline{z}}\,\Theta_{\overline{z}}\,
\rho_{ww}
+
\Theta_{\overline{z}\overline{z}}\,
\rho_w,
\notag
\\
0
&
\,\equiv\,
\rho_{\overline{z}\overline{w}}
+
\Theta_{\overline{z}}\,
\rho_{w\overline{w}}
+
\Theta_{\overline{w}}\,
\rho_{\overline{z}w}
+
\Theta_{\overline{z}}\,\Theta_{\overline{w}}\,
\rho_{ww}
+
\Theta_{\overline{z}\overline{w}}\,
\rho_w,
\notag
\\
0
&
\,\equiv\,
\rho_{\overline{w}\overline{w}}
+
2\,\Theta_{\overline{w}}\,
\rho_{w\overline{w}}
+
\Theta_{\overline{w}}\,\Theta_{\overline{w}}\,
\rho_{ww}
+
\Theta_{\overline{w}\overline{w}}\,
\rho_w.
\notag
\end{align}
It holds that:
\[
\big\{
\rho_w
\neq 0
\big\}
\,=\,
\big\{
\Theta_{\overline{w}}
\neq 0
\big\}
\eqno
{\scriptstyle{(\text{\rm in}\,M)}}.
\]

\begin{Definition}
Call $M$ {\sl spherical} if it is locally biholomorphic to:
\[
S^3
\,:=\,
\big\{
(z,w)\in\C^2
\colon\,
z\overline{z}
+
w\overline{w}
=
1
\big\}.
\]
\end{Definition}

When $M$ is connected, the principle of analytic continuation
guarantees propagation of this property.
Next, set:
\[
\Delta
\,:=\,
-\,\Theta_{\overline{w}}\,
\Theta_{z\overline{z}}
+
\Theta_{\overline{z}}\,
\Theta_{z\overline{w}}.
\]

\begin{Lemma}
At a point $p \in \{ \Theta_{\overline{w}} \neq 0\}${\em :}
\[
M\,\,
\text{is Levi nondegenerate at}\,\,
p
\,\,\,\Longleftrightarrow\,\,\,
\Delta(p)
\,\neq\,
0.\eqno\qed
\]
\end{Lemma}

Levi nondegeneracy being a biholomorphically invariant
feature, spherical $M$ are so since $S^3$ is.

Without restricting assumptions like {\em e.g.} {\sl rigidity}
or {\sl tubity} (\cite{Isaev-2011}), 
an explicit, complete characterization
of sphericity in terms of some defining function for a hypersurface
$M^3 \subset \C^2$ appeared in October 2009 as
{\footnotesize\sf arxiv.org/abs/0910.1694/},
{\em cf.} also~\cite{Nurowski-Sparling-2003, Huang-Ji-2007}.
To recall it, set:
\[
\Box
\,:=\,
\frac{\Delta}{-\,\Theta_{\overline{w}}},
\]
and use instead:
\[
\aligned
\overline{\mathcal{L}}
\,:=\,
&\,
-\,
\frac{1}{\rho_{\overline{w}}}\,
\overline{L}
\\
\,=\,
&
\frac{\partial}{\partial\overline{z}}
-
\frac{\Theta_{\overline{z}}}{\Theta_{\overline{w}}}\,
\frac{\partial}{\partial\overline{w}}.
\endaligned
\]

\begin{Theorem}
{\rm ({\cite{Merker-2010}})}
At a point $p \in \{ \Theta_{\overline{w}} \neq 0\}$, 
the hypersurface $M$ is spherical if and only if, near $p$:
\[
0
\,\equiv\,
\frac{1}{\Box}\,
\overline{\mathcal{L}}
\bigg(
\frac{1}{\Box}\,
\overline{\mathcal{L}}
\bigg(
\frac{1}{\Box}\,
\overline{\mathcal{L}}
\bigg(
\frac{1}{\Box}\,
\overline{\mathcal{L}}
\Big(
\Theta_{zz}
\Big)
\bigg)
\bigg)
\bigg).
\eqno\qed
\]
\end{Theorem}

Exchanging $z \longleftrightarrow w$ yields a similar formula at
points $q \in \{ \rho_z \neq 0\}$.

\begin{Corollary}
\label{7-terms-Theta}
In $\{ \rho_w \neq 0\} = \{ \Theta_{\overline{w}} \neq 0\}$,
a partly expanded characterization of sphericity is:
\reqnomode\usetagform{EngelLie}
\begin{align}
0
&
\,\,\equiv\,\,
\frac{\overline{\mathcal{L}}^4(\Theta_{zz})}{\Box^4}
\,-\,
\notag
\\
&
\ \ \ \ \
-\,
6\,
\frac{
\overline{\mathcal{L}}(\Box)\,
\overline{\mathcal{L}}^3(\Theta_{zz})}{
\Box^5}
-
4\,
\frac{
\overline{\mathcal{L}}^2(\Box)\,
\overline{\mathcal{L}}^2(\Theta_{zz})}{
\Box^5}
-
\frac{
\overline{\mathcal{L}}^3(\Box)\,
\overline{\mathcal{L}}(\Theta_{zz})}{
\Box^5}
\,
+
\notag
\\
&
\ \ \ \ \
+
15\,
\frac{
\big[\,\overline{\mathcal{L}}(\Box)\big]^2\,
\overline{\mathcal{L}}^2(\Theta_{zz})}{
\Box^6}
+
10\,
\frac{
\overline{\mathcal{L}}(\Box)\,
\overline{\mathcal{L}}^2(\Box)\,
\overline{\mathcal{L}}(\Theta_{zz})}{
\Box^6}
\,-
\notag
\\
&
\ \ \ \ \
-\,
15\,
\frac{
\big[\,\overline{\mathcal{L}}(\Box)\big]^3\,
\overline{\mathcal{L}}(\Theta_{zz})}{
\Box^7}.
\tag
\qed
\end{align}

\end{Corollary}

Without presenting details, it is known that Cartan's treatment
of the concerned biholomorphic equivalence problem brings
a single invariant function:
\[
\ICartan^M
\colon
\ \ \
M
\,\longrightarrow\,
\C,
\]
other invariants being (covariant) derivations of it, and that:
\[
M\,\,
\text{\rm is spherical}\,\,
\,\,\,\Longleftrightarrow\,\,\,
0
\,\equiv\,
\ICartan^M.
\]

\begin{Notation}
For two functions $\Iaux_1 \colon M \longrightarrow \C$ and 
$\Iaux_2 \colon M \longrightarrow \C$, write:
\[
\Iaux_2
\,\doteqdot\,
\Iaux_1,
\]
when there is a nowhere vanishing function $\mu \colon M 
\longrightarrow \C \backslash \{0\}$ such that:
\[
\Iaux_2
\,=\,
\mu\,
\Iaux_1.
\]
\end{Notation}

For instance:
\[
\ICartan^M
\,\doteqdot\,
\bigg(
\frac{1}{\Box}\,
\overline{\mathcal{L}}
\bigg)^{\!4}
\Big(\Theta_{zz}\Big).
\]

Now, translate the formula of Corollary~{\ref{7-terms-Theta}}
to the case where $M$ is given in implicit representation:
\[
0
\,=\,
\rho
\big(
z,w,\overline{z},\overline{w}
\big).
\]
Set:
\[
\Hessian(\rho)
\,:=\,
\rho_z\rho_z\,\rho_{ww}
-
2\,\rho_z\rho_w\,\rho_{zw}
+
\rho_w\rho_w\,\rho_{zz},
\]
with (exercise) on $\{\rho_w\neq 0\}$:
\[
\Theta_{zz}
\,=\,
-\,\frac{\Hessian(\rho)}{\rho_w\rho_w\rho_w}.
\]
Remind the {\sl Levi determinant:}
\[
\Levi(\rho)
\,:=\,
\rho_{\overline{z}}\rho_z\rho_{w\overline{w}}
-
\rho_{\overline{z}}\rho_w\rho_{z\overline{w}}
-
\rho_{\overline{w}}\rho_z\rho_{\overline{z}w}
+
\rho_{\overline{w}}\rho_w\rho_{z\overline{z}},
\]
that satisfies on $\{\rho_w \neq 0\}$:
\[
\Levi(\rho)
\,\doteqdot\,
\Delta,
\]
{\em i.e.} more precisely (exercise) thanks to~({\ref{2-rho-Theta}}):
\[
\Levi(\rho)
\,=\,
-\,
\rho_w\,\rho_w\,\rho_w\,\Delta.
\]

\begin{Corollary}
On $\{\rho_w \neq 0\}$, up to a nowhere vanishing function:
\[
\ICartan^M
\,\doteqdot\,
\Iaux_{[w]},
\]
where:
\begin{footnotesize}
\reqnomode\usetagform{EngelLie}
\begin{align}
\Iaux_{[w]}
&
\,\,:=\,\,
12\,\big(\rho_w\big)^9\,
\bigg\{
\bigg[
\frac{\Levi(\rho)}{\rho_w^2}
\bigg]^3
\overline{L}^4
\bigg(
\frac{\Hessian(\rho)}{\rho_w^3}
\bigg)
\,-
\notag
\\
&
\ \ \ \ \
-\,
6\,
\bigg[
\frac{\Levi(\rho)}{\rho_w^2}
\bigg]^2
\overline{L}
\bigg(
\frac{\Levi(\rho)}{\rho_w^2}
\bigg)
\overline{L}^3
\bigg(
\frac{\Hessian(\rho)}{\rho_w^3}
\bigg)
-
4\,
\bigg[
\frac{\Levi(\rho)}{\rho_w^2}
\bigg]^2
\overline{L}^2
\bigg(
\frac{\Levi(\rho)}{\rho_w^2}
\bigg)
\overline{L}^2
\bigg(
\frac{\Hessian(\rho)}{\rho_w^3}
\bigg)
-
\bigg[
\frac{\Levi(\rho)}{\rho_w^2}
\bigg]^2
\overline{L}^3
\bigg(
\frac{\Levi(\rho)}{\rho_w^2}
\bigg)
\overline{L}
\bigg(
\frac{\Hessian(\rho)}{\rho_w^3}
\bigg)
\,+
\notag
\\
&
\ \ \ \ \
+
15\,
\frac{\Levi(\rho)}{\rho_w^2}
\bigg[
\overline{L}
\bigg(
\frac{\Levi(\rho)}{\rho_w^2}
\bigg)
\bigg]^2
\overline{L}^2
\bigg(
\frac{\Hessian(\rho)}{\rho_w^3}
\bigg)
+
10\,
\frac{\Levi(\rho)}{\rho_w^2}
\overline{L}
\bigg(
\frac{\Levi(\rho)}{\rho_w^2}
\bigg)
\overline{L}^2
\bigg(
\frac{\Levi(\rho)}{\rho_w^2}
\bigg)
\overline{L}
\bigg(
\frac{\Hessian(\rho)}{\rho_w^3}
\bigg)
\,-
\notag
\\
&
\ \ \ \ \
-\,
15\,
\bigg[
\overline{L}
\bigg(
\frac{\Levi(\rho)}{\rho_w^2}
\bigg)
\bigg]^3
\overline{L}
\bigg(
\frac{\Hessian(\rho)}{\rho_w^3}
\bigg)
\bigg\}.
\tag{\qed}
\end{align}
\end{footnotesize}
\end{Corollary}

Furthermore, exchanging $z \longleftrightarrow w$, there is an {\em
exact} formal coincidence (exercise!):
\[
\Iaux_{[z]}
\,=\,
\Iaux_{[w]}.
\]
In~\cite{Merker-survey-2017}, an alternative formula for an equivalent
invariant $\Maux \doteqdot \Iaux_{[w]}$ is discussed, but it
incorporates $5! = 120$ terms instead of $7$ above, and is less
cleaned up or finalized to really compute exciting things (by hand!). 

%%%%%%%%%%%%%%%%%%%%%%%%%%%%%%%%%%%%%%%%%%%%%%%%%%%%%%%%%%%%%%%%%%%%%%
\Section{\bf Pullback to an Exceptional Curve on an Ellipsoid}
\label{pullback-ellipsoid}
\HEAD{\ref{pullback-ellipsoid}.~{\sf 
Pullback to an Exceptional Curve on an Ellipsoid}
}{
Wei-Guo {\sc Foo}, Jo\"el {\sc Merker}, The-Anh {\sc Ta}
Orsay University, Paris, France}
%%%%%%%%%%%%%%%%%%%%%%%%%%%%%%%%%%%%%%%%%%%%%%%%%%%%%%%%%%%%%%%%%%%%%%

To prove Theorem~{\ref{a-b-cos-sin}}, it suffices to verify that:
\[
0
\overset{\text{\bf ?}}{\,\,=\,\,}
\gamma^\ast
\big(
\Iaux_{[w]}
\big)(\theta)
\eqno
{\scriptstyle{(\forall\,\theta\,\in\,\R)}}.
\]
Drop the factor $12\, (\rho_w)^9 \doteqdot 1$, and call $\Taux_1$,
$\Taux_2$, $\Taux_3$, $\Taux_4$, $\Taux_5$, $\Taux_6$, $\Taux_7$ the
seven concerned terms, so that the goal becomes:
\[
0
\overset{\text{\bf ?}}{\,\,=\,\,}
\gamma^\ast
\big(\Taux_1\big)
+
\gamma^\ast
\big(\Taux_2\big)
+
\gamma^\ast
\big(\Taux_3\big)
+
\gamma^\ast
\big(\Taux_4\big)
+
\gamma^\ast
\big(\Taux_5\big)
+
\gamma^\ast
\big(\Taux_6\big)
+
\gamma^\ast
\big(\Taux_7\big).
\]

Hand computations provide formulas of the shape:
\[
\aligned
\Taux_1
&
\,=\,
{\textstyle{\frac{1}{8}}}\,
\isqrt\,
(a-1)\,
\frac{\Naux_1}\Daux,
\\
\Taux_2
&
\,=\,
{\textstyle{\frac{3}{4}}}\,
\isqrt\,
(a-1)\,
\frac{\Naux_2}\Daux,
\\
\Taux_3
&
\,=\,
{\textstyle{\frac{1}{2}}}\,
\isqrt\,
(a-1)\,
\frac{\Naux_3}\Daux,
\\
\Taux_4
&
\,=\,
{\textstyle{\frac{1}{8}}}\,
\isqrt\,
(a-1)\,
\frac{\Naux_4}\Daux,
\\
\Taux_5
&
\,=\,
{\textstyle{\frac{15}{8}}}\,
\isqrt\,
(a-1)\,
\frac{\Naux_5}\Daux,
\\
\Taux_6
&
\,=\,
{\textstyle{\frac{5}{4}}}\,
\isqrt\,
(a-1)\,
\frac{\Naux_6}\Daux,
\\
\Taux_7
&
\,=\,
{\textstyle{\frac{15}{8}}}\,
\isqrt\,
(a-1)\,
\frac{\Naux_7}\Daux,
\endaligned
\]
with, in denominator place: 
\[
\footnotesize
\aligned
\Daux
\,:=\,
\Big(
\sqrt{a}\,\cossmall\,\theta
-
\isqrt\,\sqrt{b}\,
\sinsmall\,\theta
\Big)^8\,
\big(
a\,b-1
\big)\,
\bigg(
\frac{
b-1}{
a\,b-1}
\bigg)^{\!\frac{11}{2}},
\endaligned
\]
with numerator 1:
\[
\!\!\!\!\!\!\!\!\!\!\!\!\!\!\!\!\!\!\!\!\!\!\!\!\!\!\!\!\!\!
\!\!\!\!\!\!\!\!\!\!
\scriptsize
\aligned
\Naux_1
:=
\cossmall^7\theta\,
&
\Big[
499\,a^{9/2}b^3
+
625\,a^{9/2}b^2
-
233\,a^{7/2}b^3
+
205\,a^{9/2}b
-
631\,a^{7/2}b^2
+
15\,a^{9/2}
-
415\,a^{7/2}b
-
65\,a^{7/2}
\Big]
\\
+\,
\isqrt\,
\cossmall^6\theta\,\sinsmall\,\theta\,
&
\Big[
2887\,a^4b^{7/2}
+
4401\,a^4b^{5/2}
-
1297\,a^3b^{7/2}
+
1905\,a^4b^{3/2}
-
4059\,a^3b^{5/2}
+
215\,a^4b^{1/2}
-
3327\,a^3b^{3/2}
-
725\,a^3b^{1/2}
\Big]
\\
+\,
\cossmall^5\theta\,\sinsmall^2\theta\,
&
\Big[
-7023\,a^{7/2}b^4
-
13021\,a^{7/2}b^3
+
3013\,a^{5/2}b^4
-
7105\,a^{7/2}b^2
+
11011\,a^{5/2}b^3
-
1075\,a^{7/2}b
+
11059\,a^{5/2}b^2
+
3141\,a^{5/2}b
\Big]
\\
+\,
\isqrt\,
\cossmall^4\theta\,\sinsmall^3\theta\,
&
\Big[
-\,9267\,a^3b^{9/2}
-
20989\,a^3b^{7/2}
+
3757\,a^2b^{9/2}
-
14101\,a^3\,b^{5/2}
+
16279\,a^2b^{7/2}
-
2683\,a^3b^{3/2}
+
19891\,a^2b^{5/2}
+
7113\,a^2b^{3/2}
\Big]
\\
+\,
\cossmall^3\theta\,\sinsmall^4\theta\,
&
\Big[
7113\,a^{5/2}b^5
+
19891\,a^{5/2}b^4
-
2683\,a^{3/2}b^5
+
16279\,a^{5/2}b^3
-
14101\,a^{3/2}b^4
+
3757\,a^{5/2}b^2
-
20989\,a^{3/2}b^3
-
9267\,a^{3/2}b^2
\Big]
\\
+\,
\isqrt\,
\cossmall^2\theta\,\sinsmall^5\theta\,
&
\Big[
3141\,a^2b^{11/2}
+
11059\,a^2b^{9/2}
-
1075\,ab^{11/2}
+
11011\,a^2b^{7/2}
-
7105\,ab^{9/2}
+
3013\,a^2b^{5/2}
-
13021\,ab^{7/2}
-
7023\,ab^{5/2}
\Big]
\\
+\,
\cossmall^1\theta\,\sinsmall^6\theta\,
&
\Big[
-\,725\,a^{3/2}b^6
-
3327\,a^{3/2}b^5
+
215\,a^{1/2}b^6
-
4059\,a^{3/2}b^4
+
1905\,a^{1/2}b^5
-
1297\,a^{3/2}b^3
+
4401\,a^{1/2}b^4
+
2287\,a^{1/2}b^3
\Big]
\\
+\,
\isqrt\,
\sinsmall^7\,\theta\,
&
\Big[
-\,65\,ab^{13/2}
-
415\,ab^{11/2}
+
15\,b^{13/2}
-
631\,ab^{9/2}
+
205\,b^{11/2}
-
233\,ab^{7/2}
+
625\,b^{9/2}
+
499\,b^{7/2}
\Big],
\endaligned
\]
with numerator 2:
\[
\!\!\!\!\!\!\!\!\!\!\!\!\!\!\!\!\!\!\!\!\!\!\!\!\!
\scriptsize
\aligned
\Naux_2
:=
\cossmall^7\theta\,
&
\Big[
-\,165\,a^{9/2}b^3
-
193\,a^{9/2}\,b^2
+
93\,a^{7/2}b^3
-
67\,a^{9/2}b
+
205\,a^{7/2}b^2
-
7\,a^{9/2}
+
115\,a^{7/2}b
+
19\,a^{7/2}
\Big]
\\
+\,
\isqrt\,
\cossmall^6\theta\,\sinsmall\,\theta\,
&
\Big[
-\,925\,a^4b^{7/2}
-
1389\,a^4b^{5/2}
+
505\,a^3b^{7/2}
-
627\,a^4b^{3/2}
+
1341\,a^3b^{5/2}
-
83\,a^4b^{1/2}
+
975\,a^3b^{3/2}
+
203\,a^3b^{1/2}
\Big]
\\
+\,
\cossmall^5\theta\,\sinsmall^2\theta\,
&
\Big[
2177\,a^{7/2}b^4
+
4141\,a^{7/2}b^3
-
1145\,a^{5/2}b^4
+
2359\,a^{7/2}b^2
-
3673\,a^{5/2}b^3
+
395\,a^{7/2}b
-
3367\,a^{5/2}b^2
-
887\,a^{5/2}b
\Big]
\\
+\,
\isqrt\,
\cossmall^4\theta\,\sinsmall^3\theta\,
&
\Big[
2777\,a^3b^{9/2}
+
6649\,a^3b^{7/2}
-
1397\,a^2b^{9/2}
+
4711\,a^3b^{5/2}
-
5449\,a^2b^{7/2}
+
983\,a^3b^{3/2}
-
6211\,a^2b^{5/2}
-
2063\,a^2b^{3/2}
\Big]
\\
+\,
\cossmall^3\theta\,\sinsmall^4\theta\,
&
\Big[
-\,2063\,a^{5/2}b^5
-
6211\,a^{5/2}b^4
+
983\,a^{3/2}b^5
-
5449\,a^{5/2}b^3
+
4711\,a^{3/2}b^4
-
1397\,a^{5/2}b^2
+
6649\,a^{3/2}b^3
+
2777\,a^{3/2}b^2
\Big]
\\
+\,
\isqrt\,
\cossmall^2\theta\,\sinsmall^5\theta\,
&
\Big[
-\,887\,a^2b^{11/2}
-
3367\,a^2b^{9/2}
+
395\,ab^{11/2}
-
3673\,a^2b^{7/2}
+
2359\,ab^{9/2}
-
1145\,a^2b^{5/2}
+
4141\,ab^{7/2}
+
2177\,ab^{5/2}
\Big]
\\
+\,
\cossmall^1\theta\,\sinsmall^6\theta\,
&
\Big[
203\,a^{3/2}b^6
+
975\,a^{3/2}b^5
-
83\,a^{1/2}b^6
+
1341\,a^{3/2}b^4
-
627\,a^{1/2}b^5
+
505\,a^{3/2}b^3
-
1389\,a^{1/2}b^4
-
925\,a^{1/2}b^3
\Big]
\\
+\,
\isqrt\,
\sinsmall^7\,\theta\,
&
\Big[
19\,ab^{13/2}
+
115\,ab^{11/2}
-
7\,b^{13/2}
+
205\,ab^{9/2}
-
67\,b^{11/2}
+
93\,ab^{7/2}
-
193\,b^{9/2}
-
165\,b^{7/2}
\Big],
\endaligned
\]
with numerator 3:
\[
\!\!\!\!\!\!\!\!\!\!\!\!\!\!\!\!\!\!\!\!
\scriptsize
\aligned
\Naux_3
:=
\cossmall^7\theta\,
&
\Big[
-\,91\,a^{9/2}b^3
-
109\,a^{9/2}b^2
+
65\,a^{7/2}b^3
-
37\,a^{9/2}b
+
115\,a^{7/2}b^2
-
3\,a^{9/2}
+
55\,a^{7/2}b
+
5\,a^{7/2}
\Big]
\\
+\,
\isqrt\,
\cossmall^6\theta\,\sinsmall\,\theta\,
&
\Big[
-\,499\,a^4b^{7/2}
-
777\,a^4b^{5/2}
+
349\,a^3b^{7/2}
-
357\,a^4b^{3/2}
+
771\,a^3b^{5/2}
-
47\,a^4b^{1/2}
+
483\,a^3b^{3/2}
+
77\,a^3b^{1/2}
\Big]
\\
+\,
\cossmall^5\theta\,\sinsmall^2\theta\,
&
\Big[
1143\,a^{7/2}b^4
+
2281\,a^{7/2}b^3
-
781\,a^{5/2}b^4
+
1369\,a^{7/2}b^2
-
2143\,a^{5/2}b^3
+
247\,a^{7/2}b
-
1723\,a^{5/2}b^2
-
393\,a^{5/2}b
\Big]
\\
+\,
\isqrt\,
\cossmall^4\theta\,\sinsmall^3\theta\,
&
\Big[
1407\,a^3b^{9/2}
+
3589\,a^3b^{7/2}
-
937\,a^2b^{9/2}
+
2761\,a^3b^{5/2}
-
3199\,a^2b^{7/2}
+
643\,a^3b^{3/2}
-
3271\,a^2b^{5/2}
-
993\,a^2b^{3/2}
\Big]
\\
+\,
\cossmall^3\theta\,\sinsmall^4\theta\,
&
\Big[
-\,993\,a^{5/2}b^5
-
3271\,a^{5/2}b^4
+
643\,a^{3/2}b^5
-
3199\,a^{5/2}b^3
+
2761\,a^{3/2}b^4
-
937\,a^{5/2}b^2
+
3589\,a^{3/2}b^3
+
1407\,a^{3/2}b^2
\Big]
\\
+\,
\isqrt\,
\cossmall^2\theta\,\sinsmall^5\theta\,
&
\Big[
-\,393\,a^2b^{11/2}
-
1723\,a^2b^{9/2}
+
247\,ab^{11/2}
-
2143\,a^2b^{7/2}
+
1369\,ab^{9/2}
-
781\,a^2b^{5/2}
+
2281\,ab^{7/2}
+
1143\,ab^{5/2}
\Big]
\\
+\,
\cossmall^1\theta\,\sinsmall^6\theta\,
&
\Big[
77\,a^{3/2}b^6
+
483\,a^{3/2}b^5
-
47\,a^{1/2}b^6
+
771\,a^{3/2}b^4
-
357\,a^{1/2}b^5
+
349\,a^{3/2}b^3
-
777\,a^{1/2}b^4
-
499\,a^{1/2}b^3
\Big]
\\
+\,
\isqrt\,
\sinsmall^7\,\theta\,
&
\Big[
5\,ab^{13/2}
+
55\,ab^{11/2}
-
3\,b^{13/2}
+
115\,ab^{9/2}
-
37\,b^{11/2}
+
65\,ab^{7/2}
-
109\,b^{9/2}
-
91\,b^{7/2}
\Big],
\endaligned
\]
with numerator 4:
\[
\!\!\!\!\!\!\!\!\!\!\!\!\!\!\!\!\!\!\!\!
\scriptsize
\aligned
\Naux_4
:=
\cossmall^7\theta\,
&
\Big[
-\,75\,a^{9/2}b^3
-
91\,a^{9/2}b^2
+
75\,a^{7/2}b^3
-
25\,a^{9/2}b
+
91\,a^{7/2}b^2
-
a^{9/2}
+
25\,a^{7/2}b
+
a^{7/2}
\Big]
\\
+\,
\isqrt\,
\cossmall^6\theta\,\sinsmall\,\theta\,
&
\Big[
-\,391\,a^4b^{7/2}
-
639\,a^4b^{5/2}
+
391\,a^3b^{7/2}
-
285\,a^4b^{3/2}
+
639\,a^3b^{5/2}
-
29\,a^4b^{1/2}
+
285\,a^3b^{3/2}
+
29\,a^3b^{1/2}
\Big]
\\
+\,
\cossmall^5\theta\,\sinsmall^2\theta\,
&
\Big[
839\,a^{7/2}b^4
+
1831\,a^{7/2}b^3
-
839\,a^{5/2}b^4
+
1165\,a^{7/2}b^2
-
1831\,a^{5/2}b^3
+
197\,a^{7/2}b
-
1165\,a^{5/2}b^2
-
197\,a^{5/2}b
\Big]
\\
+\,
\isqrt\,
\cossmall^4\theta\,\sinsmall^3\theta\,
&
\Big[
947\,a^3b^{9/2}
+
2779\,a^3b^{7/2}
-
947\,a^2b^{9/2}
+
2401\,a^3b^{5/2}
-
2779\,a^2b^{7/2}
+
593\,a^3b^{3/2}
-
2401\,a^2b^{5/2}
-
593\,a^2b^{3/2}
\Big]
\\
+\,
\cossmall^3\theta\,\sinsmall^4\theta\,
&
\Big[
-\,593\,a^{5/2}b^5
-
2401\,a^{5/2}b^4
+
593\,a^{3/2}b^5
-
2779\,a^{5/2}b^3
+
2401\,a^{3/2}b^4
-
947\,a^{5/2}b^2
+
2779\,a^{3/2}b^3
+
947\,a^{3/2}b^2
\Big]
\\
+\,
\isqrt\,
\cossmall^2\theta\,\sinsmall^5\theta\,
&
\Big[
-\,197\,a^2b^{11/2}
-
1165\,a^2b^{9/2}
+
197\,ab^{11/2}
-
1831\,a^2b^{7/2}
+
1165\,ab^{9/2}
-
839\,a^2b^{5/2}
+
1831\,ab^{7/2}
+
839\,ab^{5/2}
\Big]
\\
+\,
\cossmall^1\theta\,\sinsmall^6\theta\,
&
\Big[
29\,a^{3/2}b^6
+
285\,a^{3/2}b^5
-
29\,a^{1/2}b^6
+
639\,a^{3/2}b^4
-
285\,a^{1/2}b^5
+
391\,a^{3/2}b^3
-
639\,a^{1/2}b^4
-
391\,a^{1/2}b^3
\Big]
\\
+\,
\isqrt\,
\sinsmall^7\,\theta\,
&
\Big[
ab^{13/2}
+
25\,ab^{11/2}
-
b^{13/2}
+
91\,ab^{9/2}
-
25\,b^{11/2}
+
75\,ab^{7/2}
-
91\,b^{9/2}
-
75\,b^{7/2}
\Big],
\endaligned
\]
with numerator 5:
\[
\!\!\!\!\!\!\!\!\!\!\!\!\!\!\!\!\!\!\!\!
\scriptsize
\aligned
\Naux_5
:=
\cossmall^7\theta\,
&
\Big[
63\,a^{9/2}b^3
+
69\,a^{9/2}b^2
-
45\,a^{7/2}b^3
+
25\,a^{9/2}b
-
75\,a^{7/2}b^2
+
3\,a^{9/2}
-
35\,a^{7/2}b
-
5\,a^{7/2}
\Big]
\\
+\,
\isqrt\,
\cossmall^6\theta\,\sinsmall\,\theta\,
&
\Big[
339\,a^4b^{7/2}
+
509\,a^4b^{5/2}
-
237\,a^3b^{7/2}
+
237\,a^4b^{3/2}
-
511\,a^3b^{5/2}
+
35\,a^4b^{1/2}
-
315\,a^3b^{3/2}
-
57\,a^3b^{1/2}
\Big]
\\
+\,
\cossmall^5\theta\,\sinsmall^2\theta\,
&
\Big[
-\,763\,a^{7/2}b^4
-
1521\,a^{7/2}b^3
+
521\,a^{5/2}b^4
-
909\,a^{7/2}b^2
+
1431\,a^{5/2}b^3
-
167\,a^{7/2}b
+
1143\,a^{5/2}b^2
+
265\,a^{5/2}b
\Big]
\\
+\,
\isqrt\,
\cossmall^4\theta\,\sinsmall^3\theta\,
&
\Big[
-\,927\,a^3b^{9/2}
-
2409\,a^3b^{7/2}
+
617\,a^2b^{9/2}
-
1841\,a^3b^{5/2}
+
2139\,a^2b^{7/2}
-
423\,a^3b^{3/2}
+
2191\,a^2b^{5/2}
+
653\,a^2b^{3/2}
\Big]
\\
+\,
\cossmall^3\theta\,\sinsmall^4\theta\,
&
\Big[
653\,a^{5/2}b^5
+
2191\,a^{5/2}b^4
-
423\,a^{3/2}b^5
+
2139\,a^{5/2}b^3
-
1841\,a^{3/2}b^4
+
617\,a^{5/2}b^2
-
2409\,a^{3/2}b^3
-
927\,a^{3/2}b^2
\Big]
\\
+\,
\isqrt\,
\cossmall^2\theta\,\sinsmall^5\theta\,
&
\Big[
265\,a^2b^{11/2}
+
1143\,a^2b^{9/2}
-
167\,ab^{11/2}
+
1431\,a^2b^{7/2}
-
909\,ab^{9/2}
+
521\,a^2b^{5/2}
-
1521\,ab^{7/2}
-
763\,ab^{5/2}
\Big]
\\
+\,
\cossmall^1\theta\,\sinsmall^6\theta\,
&
\Big[
-\,57\,a^{3/2}b^6
-
315\,a^{3/2}b^5
+
35\,a^{1/2}b^6
-
511\,a^{3/2}b^4
+
237\,a^{1/2}b^5
-
237\,a^{3/2}b^3
+
509\,a^{1/2}b^4
+
339\,a^{1/2}b^3
\Big]
\\
+\,
\isqrt\,
\sinsmall^7\,\theta\,
&
\Big[
-\,5\,ab^{13/2}
-
35\,ab^{11/2}
+
3\,b^{13/2}
-
75\,ab^{9/2}
+
25\,b^{11/2}
-
45\,ab^{7/2}
+
69\,b^{9/2}
+
63\,b^{7/2}
\Big],
\endaligned
\]
with numerator 6:
\[
\!\!\!\!\!\!\!\!\!\!\!\!\!\!\!\!\!\!\!\!
\scriptsize
\aligned
\Naux_6
:=
\cossmall^7\theta\,
&
\Big[
39\,a^{9/2}b^3
+
43\,a^{9/2}b^2
-
39\,a^{7/2}b^3
+
13\,a^{9/2}b
-
43\,a^{7/2}b^2
+
a^{9/2}
-
13\,a^{7/2}b
-
a^{7/2}
\Big]
\\
+\,
\isqrt\,
\cossmall^6\theta\,\sinsmall\,\theta\,
&
\Big[
199\,a^4b^{7/2}
+
315\,a^4b^{5/2}
-
199\,a^3b^{7/2}
+
141\,a^4b^{3/2}
-
315\,a^3b^{5/2}
+
17\,a^4b^{1/2}
-
141\,a^3b^{3/2}
-
17\,a^3b^{1/2}
\Big]
\\
+\,
\cossmall^5\theta\,\sinsmall^2\theta\,
&
\Big[
-\,419\,a^{7/2}b^4
-
919\,a^{7/2}b^3
+
419\,a^{5/2}b^4
-
577\,a^{7/2}b^2
+
919\,a^{5/2}b^3
-
101\,a^{7/2}b
+
577\,a^{5/2}b^2
+
101\,a^{5/2}b
\Big]
\\
+\,
\isqrt\,
\cossmall^4\theta\,\sinsmall^3\theta\,
&
\Big[
-\,467\,a^3b^{9/2}
-
1399\,a^3b^{7/2}
+
467\,a^2b^{9/2}
-
1201\,a^3b^{5/2}
+
1399\,a^2b^{7/2}
-
293\,a^3b^{3/2}
+
1201\,a^2b^{5/2}
+
293\,a^2b^{3/2}
\Big]
\\
+\,
\cossmall^3\theta\,\sinsmall^4\theta\,
&
\Big[
293\,a^{5/2}b^5
+
1201\,a^{5/2}b^4
-
293\,a^{3/2}b^5
+
1399\,a^{5/2}b^3
-
1201\,a^{3/2}b^4
+
467\,a^{5/2}b^2
-
1399\,a^{3/2}b^3
-
467\,a^{3/2}b^2
\Big]
\\
+\,
\isqrt\,
\cossmall^2\theta\,\sinsmall^5\theta\,
&
\Big[
101\,
a^2b^{11/2}
+
577\,a^2b^{9/2}
-
101\,ab^{11/2}
+
919\,a^2b^{7/2}
-
577\,ab^{9/2}
+
419\,a^2b^{5/2}
-
919\,ab^{7/2}
-
419\,ab^{5/2}
\Big]
\\
+\,
\cossmall^1\theta\,\sinsmall^6\theta\,
&
\Big[
-\,17\,a^{3/2}b^6
-
141\,a^{3/2}b^5
+
17\,a^{1/2}b^6
-
315\,a^{3/2}b^4
+
141\,a^{1/2}b^5
-
199\,a^{3/2}b^3
+
315\,a^{1/2}b^4
+
199\,a^{1/2}b^3
\Big]
\\
+\,
\isqrt\,
\sinsmall^7\,\theta\,
&
\Big[
-\,ab^{13/2}
-
13\,ab^{11/2}
+
b^{13/2}
-
43\,ab^{9/2}
+
13\,b^{11/2}
-
39\,ab^{7/2}
+
43\,b^{9/2}
+
39\,b^{7/2}
\Big],
\endaligned
\]
with numerator 7:
\[
\!\!\!\!\!\!\!\!\!\!\!\!\!\!\!\!\!\!\!\!
\scriptsize
\aligned
\Naux_7
:=
\cossmall^7\theta\,
&
\Big[
-\,27\,a^{9/2}b^3
-
27\,a^{9/2}b^2
+
27\,a^{7/2}b^3
-
9\,a^{9/2}b
+
27\,a^{7/2}b^2
-
a^{9/2}
+
9\,a^{7/2}b
+
a^{7/2}
\Big]
\\
+\,
\isqrt\,
\cossmall^6\theta\,\sinsmall\,\theta\,
&
\Big[
-\,135\,a^4b^{7/2}
-
207\,a^4b^{5/2}
+
135\,a^3b^{7/2}
-
93\,a^4b^{3/2}
+
207\,a^3b^{5/2}
-
13\,a^4b^{1/2}
+
93\,a^3b^{3/2}
+
13\,a^3b^{1/2}
\Big]
\\
+\,
\cossmall^5\theta\,\sinsmall^2\theta\,
&
\Big[
279\,a^{7/2}b^4
+
615\,a^{7/2}b^3
-
279\,a^{5/2}b^4
+
381\,a^{7/2}b^2
-
615\,a^{5/2}b^3
+
69\,a^{7/2}b
-
381\,a^{5/2}b^2
-
69\,a^{5/2}b
\Big]
\\
+\,
\isqrt\,
\cossmall^4\theta\,\sinsmall^3\theta\,
&
\Big[
307\,a^3b^{9/2}
+
939\,a^3b^{7/2}
-
307\,a^2b^{9/2}
+
801\,a^3b^{5/2}
-
939\,a^2b^{7/2}
+
193\,a^3b^{3/2}
-
801\,a^2b^{5/2}
-
193\,a^2b^{3/2}
\Big]
\\
+\,
\cossmall^3\theta\,\sinsmall^4\theta\,
&
\Big[
-\,193\,a^{5/2}b^5
-
801\,a^{5/2}b^4
+
193\,a^{3/2}b^5
-
939\,a^{5/2}b^3
+
801\,a^{3/2}b^4
-
307\,a^{5/2}b^2
+
939\,a^{3/2}b^3
+
307\,a^{3/2}b^2
\Big]
\\
+\,
\isqrt\,
\cossmall^2\theta\,\sinsmall^5\theta\,
&
\Big[
-\,69\,a^2b^{11/2}
-
381\,a^2b^{9/2}
+
69\,ab^{11/2}
-
615\,a^2b^{7/2}
+
381\,ab^{9/2}
-
279\,a^2b^{5/2}
+
615\,ab^{7/2}
+
279\,ab^{5/2}
\Big]
\\
+\,
\cossmall^1\theta\,\sinsmall^6\theta\,
&
\Big[
13\,a^{3/2}b^6
+
93\,a^{3/2}b^5
-
13\,a^{1/2}b^6
+
207\,a^{3/2}b^4
-
93\,a^{1/2}b^5
+
135\,a^{3/2}b^3
-
207\,a^{1/2}b^4
-
135\,a^{1/2}b^3
\Big]
\\
+\,
\isqrt\,
\sinsmall^7\,\theta\,
&
\Big[
ab^{13/2}
+
9\,ab^{11/2}
-
b^{13/2}
+
27\,ab^{9/2}
-
9\,b^{11/2}
+
27\,ab^{7/2}
-
27\,b^{9/2}
-
27\,b^{7/2}
\Big].
\endaligned
\]

\proof[End of proof of Theorem~{\ref{a-b-cos-sin}}]
The sum:
\[
{\textstyle{\frac{1}{8}}}\,
\Naux_1
(\theta)
+
{\textstyle{\frac{3}{4}}}\,
\Naux_2
(\theta)
+
{\textstyle{\frac{1}{2}}}\,
\Naux_3
(\theta)
+
{\textstyle{\frac{1}{8}}}\,
\Naux_4
(\theta)
+
{\textstyle{\frac{15}{8}}}\,
\Naux_5
(\theta)
+
{\textstyle{\frac{5}{4}}}\,
\Naux_6
(\theta)
+
{\textstyle{\frac{15}{8}}}\,
\Naux_7
(\theta)
\,=\,
0,
\]
is indeed (visually!) identically null.
\endproof

%%%%%%%%%%%%%%%%%%%%%%%%%%%%%%%%%%%%%%%%%%%%%%%%%%%%%%%%%%%%%%%%%%%%%%

\newpage

%%%%%%%%%%%%%%%%%%%%%%%%%%%%%%%%%%%%%%%%%%%%%%%%%%%%%%%%%%%%%%%%%%%%%%
%%%%%%%%%%%%%%%%%%%%%%%%%%%%%%%%%%%%%%%%%%%%%%%%%%%%%%%%%%%%%%%%%%%%%%
%%%%%%%%%%%%%%%%%%%%%%%%%%%%%%%%%%%%%%%%%%%%%%%%%%%%%%%%%%%%%%%%%%%%%%

\bibliographystyle{amsplain}

\noindent
{\scriptsize\sf fooweiguo@gmail.com, 
joel.merker@math.u-psud.fr, the-anh.ta@math.u-psud.fr}

%%%%%%%%%%%%%%%%%%%%%%%%%%%%%%%%%%%%%%%%%%%%%%%%%%%%%%%%%%%%%%%%%%%%%%

\vfill\end{document}

%% file: macros.tex
\newtheorem{Theorem}[equation]{Theorem}

\newtheorem{Lemma}[equation]{Lemma}

\newtheorem{Corollary}[equation]{Corollary}

\theoremstyle{definition}

\newtheorem{Definition}[equation]{Definition}
\newtheorem{Notation}[equation]{Notation}

\newtheorem{Question}[equation]{Question}

\newcommand{\C}{\mathbb{C}}

\newcommand{\R}{\mathbb{R}}

\newcommand{\NN}{\text{\sc n}}

\newcommand{\xaux}{{\text{\usefont{T1}{qcs}{m}{sl}x}}}
\newcommand{\yaux}{{\text{\usefont{T1}{qcs}{m}{sl}y}}}
\newcommand{\zaux}{{\text{\usefont{T1}{qcs}{m}{sl}z}}}

\newcommand{\Daux}{{\text{\usefont{T1}{qcs}{m}{sl}D}}}

\newcommand{\Iaux}{{\text{\usefont{T1}{qcs}{m}{sl}I}}}

\newcommand{\Maux}{{\text{\usefont{T1}{qcs}{m}{sl}M}}}
\newcommand{\Naux}{{\text{\usefont{T1}{qcs}{m}{sl}N}}}
\newcommand{\Taux}{{\text{\usefont{T1}{qcs}{m}{sl}T}}}

\definecolor{blue}{cmyk}{1.,1.,0.,0.63}
\definecolor{red}{cmyk}{0.,1.,1.,0.63}
\definecolor{green}{cmyk}{1.,0.,1.,0.63}
\definecolor{black}{cmyk}{1.,1.,1.,1.}

\newcommand{\blue}{\textcolor{blue}}

\makeatletter
\renewcommand{\@fnsymbol}[1]
{\ensuremath{\ifcase#1\or $*$\or $**$\or $***$\or $****$\or $*****$
\else\@ctrerr\fi}}
\makeatother

\newcommand{\HEAD}[2]{%
\pagestyle{fancy}
\fancyhead[RO]{\tiny\sf\thepage}
\fancyhead[CO]{{\tiny\sf #1}}
\fancyhead[LE]{\tiny\sf\thepage}
\fancyhead[CE]{{\tiny\sf #2}}
\fancyfoot{}}

\newcommand{\CITATION}[1]{\smallskip\hfill
\begin{minipage}[t]{13.5cm}\baselineskip=0.43cm\parindent=0.91cm
\blue{\footnotesize{\sf \!\!\!\!\!\!\!#1}}\end{minipage}\medskip}

\numberwithin{equation}{section}

\newcommand{\Section}[1]{
\renewcommand{\thesection}{\bf\arabic{section}}
\section{#1}
\renewcommand{\thesection}{\arabic{section}}}

\newcommand{\style}[1]{\text{\footnotesize{\sf #1}}}

\newcommand{\stylesmall}[1]{{\sf #1}}

\renewcommand{\cos}{\style{cos}}
\newcommand{\cossmall}{\stylesmall{cos}}

\renewcommand{\dim}{\style{dim}}

\newcommand{\Hessian}{\style{H}}

\renewcommand{\Im}{\style{Im}}

\newcommand{\Levi}{\style{L}}

\renewcommand{\lim}{\style{lim}}

\renewcommand{\Re}{\style{Re}}

\renewcommand{\sin}{\style{sin}}
\newcommand{\sinsmall}{\stylesmall{sin}}

\newcommand{\ICartan}{\mathfrak{I}_{\sf Cartan}}

\newcommand{\isqrt}{{\scriptstyle{\sqrt{-1}}}}

\newcommand{\UmbCR}{{\sf Umb}_{\sf CR}}

\newcommand{\vf}{\vfill\end{document}}